\documentclass[a4paper,11pt]{article}
\usepackage{amsmath,amssymb,amsthm,bm,bbm}
\usepackage{amscd}
\usepackage[monochrome]{color}

\usepackage{comment}

\renewcommand{\epsilon}{\varepsilon}
\renewcommand{\phi}{\varphi}

\bmdefine{\be}{e}
\bmdefine{\bE}{E}
\bmdefine{\bh}{h}
\bmdefine{\bg}{g}
\bmdefine{\bk}{k}
\bmdefine{\bp}{p}
\bmdefine{\bx}{x}
\bmdefine{\by}{y}
\bmdefine{\bs}{s}
\bmdefine{\bt}{t}
\bmdefine{\bu}{u}
\bmdefine{\bv}{v}
\bmdefine{\bw}{w}
\bmdefine{\bep}{\epsilon}

\newcommand{\B}{{\mathcal B}}
\newcommand{\R}{{\mathbb R}}

\newcommand{\Q}{{\mathbb Q}}
\newcommand{\Z}{{\mathbb Z}}
\newcommand{\C}{{\mathbb C}}
\newcommand{\N}{{\mathbb N}}
\newcommand{\LL}{{\mathcal L}}

\newcommand{\A}{{\mathcal A}}

\renewcommand{\Re}{{\mathrm{Re}}}

\newcommand{\Res}{{\mathrm{Res}}}

\newcommand{\1}{\mathbbm{1}}

\newcommand{\ov}{\overline}

\newtheorem{dfn}{DEFINITION}
\newtheorem{theorem}{THEOREM}
\newtheorem{lem}{LEMMA}

\newtheorem{cor}{COROLLARY}
\newtheorem{prop}{PROPOSITION}
\newtheorem{rem}{REMARK}

\newcommand{\Romannum}[1]{\uppercase\expandafter{\romannumeral#1}}
\makeatletter
  
  \@addtoreset{equation}{section}
\makeatother

\begin{document}
\title{Exponential Diophantine approximation and symbolic dynamics}
\author{Shigeki Akiyama, Teturo Kamae, Hajime Kaneko}
\date{}
\maketitle

\begin{abstract}
We extend the key formula in \cite{Ak-Ka:21, kan1, kan2} which intertwines 
multiplicative Markoff-Lagrange spectrum and symbolic dynamics.  
The proof uses complex analysis and elucidates the strategy of the problem. 
Moreover, the new method applies to a wide variety of polynomials possibly having multiple roots. 
We derive several consequences of this formula, which are 
expected for the Markoff-Lagrange spectrum.  
\end{abstract}

\section{Introduction}\label{section1}
For any real number $x$, we denote its integral and fractional parts by 
$\lfloor x\rfloor$ and $\{x\}$, respectively. 
Let $\|x\|$ be the distance from $x$ to the nearest integer. There exist a unique integer $u(x)$ and 
a unique real number $\epsilon(x)\in[-1/2,1/2)$ 
with $x=u(x)+\epsilon(x)$. Note that 
$\|x\|=|\epsilon(x)|$. 
For $\alpha>1$, we wish to study the set  
\[
\LL(\alpha):=
\left\{ \left. \limsup_{n\to\infty}  \| \xi\alpha^n\| \
\right| \ \xi\in \R\right\}.
\]
This type of exponential Diophantine problem seems to be beyond our reach for a general $\alpha$. 
Regardless of this situation, 
in the previous paper \cite{Ak-Ka:21} we showed 
that for a certain Pisot 
number $\alpha$ that $\LL(\alpha)$ has analogous 
property 
as the Lagrange spectrum.
The key to the proofs of \cite{Ak-Ka:21} was the formula that lift the problem to a
symbolic dynamical setting. Here is a schematic view:
\[
  \begin{CD}
  \Omega  @>{\sigma}>>  \Omega\\
  @V{h}VV    @V{h}VV \\
             \xi \alpha^n \pmod{\Z}@>{\times \alpha}>> \xi \alpha^{n+1} \pmod{\Z}
  \end{CD}
\]
where $\Omega$ is a shift space and $\sigma((t_m))=(t_{m+1})$ 
is the shift map acting on bi-infinite sequences $(t_m)_{m\in \Z}$
over a finite number of symbols. The function $h$ intertwines our problem with the one in a different setting, and
our key formula gives the commutativity of this diagram.  Although we write $\times \alpha$ symbolically 
as a function, we should put stress on the fact that 
when $\alpha\not \in \Z$, there is no appropriate compact metric space nor probability space where the multiplication $\times \alpha$ naturally acts,
which is our main difficulty to treat this exponential Diophantine problem. The direct proof of the formula was given in \cite{Ak-Ka:21}
whose underlining principle has been developed in \cite{kan1, kan2}. 
It is of intimate interest to treat $\LL(\alpha)$ for general algebraic numbers. Therefore we are 
interested in giving a unified and simpler proof of the formula which would apply to the general situation.

In this paper, we present such an expected generalization 
which is understood as an inverse matrix in infinite-dimensional vector space. 
As a byproduct
we can discuss a pretty general linear recurrence sequence whose characteristic polynomial may even have multiple roots. We shall discuss 
$$
\left\{
\limsup_{n\to \infty} \Vert  \Re (x_n)  \Vert \right\}
$$
where $x_n$ extends over complex sequences
generated by a fixed linear recurrence and give results that generalize \cite{Ak-Ka:21}.
Because of the general feature of our treatise, we are satisfied by a little weaker results than those in \cite{Ak-Ka:21}. 

We finish this introduction with two concrete examples which illustrate our results. 
Let $$f(X)=X^3 +2 X^2 +6 X - 2$$
and $k$ be a non negative integer.
Then $f(X)$ has two complex roots $\alpha, \overline{\alpha}$ and a real root $\beta$ with $\alpha\approx-1.1495 + 2.3165  \sqrt{-1}$ and $\beta\approx 0.2991$.
The set 
$$
\LL_1=\left\{\left. \limsup_{n \to \infty} \left\Vert \Re (p(n)\alpha^n) \right\Vert\ \right|\ p(x)\in \C[x],\ \mathrm{deg}(p)\le k \right\}
$$
is closed by Theorem \ref{thm:2-1}. By Theorem \ref{thm:iso}, the point
$0$ is an isolated point. One can confirm (\ref{eqn:3-19}), and thus there exists proper intervals in $\LL_1$ whose end points converge to $1/2$
by Theorem \ref{thm:pos_int_gen}. 
Here an interval is proper if it has a positive length.\footnote{Detailed discussion allows 
us to find a proper 
interval $[v,1/2]\subset \LL_1$, see Remark \ref{Interval}.} 
Though $\alpha$ is not a Pisot number nor a unit, we observe a multiplicative analogy
of the Markoff-Lagrange spectrum. \par
Similarly for a polynomial $X^2-20X+82$, the set 
$$
\LL_2=\left\{\left. \limsup_{n \to \infty} \left\Vert \xi_1(10+3\sqrt{2})^n +\xi_2(10-3\sqrt{2})^n \right\Vert\ \right|\ \xi_1,\xi_2 \in \R \right\}
$$
is closed and $0$ is an isolated point. 
Moreover, Corollary \ref{cor:disc2} gives a very 
precise description of the discrete part of the spectrum below the smallest accumulation point. 
Also one can confirm (\ref{eqn:3-19}) and there are proper intervals in $\LL_2$ whose endpoints converge to $1/2$.

\section{Basic notation}\label{section2}

For any complex number $c$, we denote its complex conjugate by $\ov{c}$. 
Moreover, for any polynomial $h(X)=\sum_{j=0}^{\ell} c_j X^j\in \C[X]$ in one variable $X$, 
let $\ov{h(X)}=\sum_{j=0}^{\ell} \ov{c_j} X^j$. \par
We introduce notation which we use in the rest of this paper. 
Let $P(X)=a_d X^d+a_{d-1}X^{d-1}+\cdots+a_0\in \Z[X]$ be 
\textcolor{red}{a polynomial} with 
$d\geq 1$ and $a_d\geq 1$. Assume that $P(X)$ satisfies the following 
three assumptions: 
\begin{enumerate}
\item $\mathrm{gcd}(a_d,a_{d-1},\ldots,a_0)=1$ and $a_0\ne 0$. 
\item $P(X)$ does not have multiple roots. 
\item {Any irreducible factor $Q(X)\in \Z[X]$ of 
$P(X)$ has at least one root 
with modulus greater than 1. }
\end{enumerate}
Moreover, we suppose that there exist $p=r_1+2r_2$ expanding roots 
and \textcolor{red}{$d-p=r'_1+2r'_2$} non expanding roots 
among the set of roots of \textcolor{red}{$P(X)=0$}, 
$\alpha_1,\ldots,\alpha_d$, i.e.,
$$|\alpha_i|>1 \text{ for } i=1,\dots, p \text{ and } |\alpha_i|\le 1 \text{ for } i=p+1,\dots d,$$
$$
\alpha_i\in \R \text { for } i=1,2,\dots, r_1,
$$
$$
\alpha_{r_1+i}\in \C\setminus \R \text{ and } \alpha_{r_1+r_2+i}=\overline{\alpha_{i+r_1}}
\text{ for } i=1,2,\dots, r_2,
$$
$$
\alpha_{p+i}\in \R \text{ for } i=1,2,\dots, r'_1,
$$
$$
\alpha_{p+r'_1+i}\in \C\setminus \R \text{ and } \alpha_{p+r'_1+r'_2+i}=\overline{\alpha_{p+r'_1+i}}
\text{ for } i=1,2,\dots, r'_2.
$$
For a nonnegative integer $k$, we let $\Xi_k$ be the set of 
\[
\bg=(g_1(X),\ldots,g_{r_1}(X),g_{1+r_1}(X),\ldots,g_p(X))\in \R[X]^{r_1}\times \C[X]^{2r_2},\] where 
$g_1(X),\ldots,g_p(X)$ are polynomials in one variable $X$ with $\deg g_i(X)\leq~k$ 
for any $i=1,\ldots,p$ and $\ov{g_{i+r_1}(X)}=g_{i+r_1+r_2}(X)$ for any 
$i=1,\ldots,r_2$. 
In the case of $k=0$, we put 
$\Xi:=\Xi_0\subset \R^{r_1}\times \C^{2r_2}$. \par
For any $\bg\in \Xi_k$, we define the linear recurrence 
$\bx=\bx(\bg)=(x_n(\bg))_{n\in\Z}$ by 
\begin{align}\label{basicdef}
x_n(\bg)=\sum_{i=1}^{p} g_i(n) \alpha_i^n\in \R.
\end{align}
\begin{lem}\label{lem:2-1}
Let $(g_1(X),\ldots,g_p(X))\in \C[X]^{r_1}\times \C[X]^{2r_2}$. 
Assume that $\deg g_j(X)\leq k$ for any $1\leq j\leq p$. \\
(1) 
\[
\sum_{j=1}^p g_j(n)\alpha_j^n=0 \mbox{ for any }
n\in \Z \mbox{ with }0\leq n\leq -1+(k+1)d
\] 
if and only if $(g_1(X),\ldots,g_p(X))=(0,\ldots,0)$. \\
(2)
\[
\sum_{j=1}^p g_j(n)\alpha_j^n\in \R \mbox{ for any }n\in \Z
\]
if and only if $(g_1(X),\ldots,g_p(X))\in \Xi_k$. 
\end{lem}
\begin{proof}
The second statement of Lemma \ref{lem:2-1} follows from 
the first statement. We now check the first statement. 
We define the square matrix $W$ of order $(k+1)d$ by 
\[
W=(W_1 W_2 \ldots W_d),
\]
where $W_{\ell}$ ($1\leq \ell\leq d$) is the $(k+1)d \times (k+1)$ matrix 
defined by 
\[
W_{\ell}=(m^n\alpha_{\ell}^m)_{0\leq m\leq -1+(k+1)d, 0\leq n\leq k}. 
\]
Using the formula for the generalized Vandermonde determinants 
proved by Flowe and Harris 
\cite{flo}, we see that $\det W\ne 0$, which implies the first statement. 
\end{proof}

\section{Intertwining formula to symbolic dynamics}\label{section3}

In this section, we describe the main tool of this paper. 

Consider the set of the sequences in (\ref{basicdef}) 
for any $\bg\in\Xi_k$. 
We consider them in modulo $\Z$, that is, as $(\epsilon_m)_{m\in \Z}=(\epsilon(x_m(\bg)))_{m\in \Z}$.
Denote 
\begin{align}\label{eqn:3-01}
\bE_f=\{(\epsilon(x_m(\bg)))_{m\in \Z};~\bg\in\Xi_k\}. 
\end{align}
The closure of $\bE_f$ define a dynamical system together with the shift: 
$$(\epsilon_m)_{m\in \Z}\mapsto(\epsilon_{m+1})_{m\in \Z}.$$
This system has a symbolic representation given by (\ref{eqn:3-2}). 
Our Theorem \ref{thm:3-1} reads 
there exist a bi-infinite sequence
$\rho_n^{f}$ ($n\in \Z$) which vanishes exponentially in the direction $n\to\infty$ 
and is at most of polynomial order of $n$ in the direction $n\to-\infty$, 
determined only by $f(X)$, that gives the inverse of this symbolic representation. 
That is, if the symbolic representation of $(\epsilon_m)$ is $(s_m)$, then
%
\begin{align}\label{eqn:3-17.5}
\epsilon_m=\sum_{j\in \Z}\rho_{-j}^f s_{j+m}
\end{align}
holds for any $m\in \Z$.
Moreover, when $P(X)$ is monic, let 
\begin{align}\label{3-0}
&\Omega_0\mbox{ be the set of sequences $(t_m)_{m\in\Z}$ of 
integers such that}\nonumber\\
&\mbox{there exists }M<0\mbox{ satisfying }t_{m}=0
\mbox{ for any }m<M,\\
&\mbox{and }t_m=O(\lambda^m)\mbox{ as }m\to\infty\mbox{ for any }
\lambda>1. \nonumber
\end{align}
Then for any $(t_m)_{m\in\Z}\in\Omega_0$, there exists 
$(x_m(\bg))_{m\in\Z}$ such that 
\begin{align}\label{eqn:3-18}
\epsilon(x_m(\bg))
\equiv
\sum_{j\in \Z}\rho_{-j}^f t_{j+m}
\pmod{\Z}
\end{align}
for any $m\in \Z$, which is our ``inverse formula" together with (\ref{eqn:3-17.5}). 

\begin{figure}[ht]
\setlength{\unitlength}{1mm}
\begin{picture}(60,50)(-20,-25)
\put(20,0){\vector(1,0){40}}
\put(60,-2){\vector(-1,0){40}}

\put(11,16){\vector(0,-1){14}}
\put(13,2){\vector(0,1){14}}

\put(19,19){\vector(1,0){40}}
\put(63,16){\vector(0,-1){14}}

\put(37,15){$\B_f$}
\put(-14,8){$\exists\pi$:factor map}
\put(15,8){inclusion map}
\put(61,18){$\R^\Z$}
\put(65,8){mod $\Z$}

\put(7,-2){$\{(s_m)\}$}
\put(62,-2){$\bE_f$}
\put(33,2){$\B_f=\A_f^{-1}$}
\put(33,-6){$\A_f=\B_f^{-1}$}
\put(10,18){$\Omega_0$}
\end{picture}
\vspace{-5em}
\caption{All the mappings commute with the shift.}
\end{figure}
\par

\subsection{Formulae by matrix with the index set $\Z\times \Z$}\label{subsection3-1}
For any nonnegative integer $k$, let 
\[
P(X)^{1+k}=:f(X)=A_DX^D+A_{D-1}X^{D-1}+\cdots+A_0\in \Z[X],
\]
where $D=(k+1)d$. 
It is known for any $\bg\in \Xi_k$ that 
$(x_m(\bg))_{n\in \Z}$ satisfies 
\begin{align}\label{eqn:3-1}
A_0 x_m(\bg)+A_1 x_{m+1}(\bg)+\cdots+A_D x_{m+D}(\bg)=0
\end{align}
for any integer $m$. 
Let 
\begin{equation}\label{eqn:3-0}
u_m=u_m(\bg):=u(x_m(\bg))\in\Z,\quad
\epsilon_m=\epsilon_m(\bg):=\epsilon(x_m(\bg)), 
\end{equation}
where for $x\in\R$, we denote 
\begin{equation}\label{eqn3-00}
x=u(x)+\epsilon(x)\mbox{ with }u(x)\in\Z\mbox{ and }\epsilon(x)\in [1/2,-1/2). 
\end{equation}
Also, let
\[
s_m=s_m(\bg):=A_0 u_m+A_1 u_{m+1}+\cdots+A_D u_{m+D}. 
\]
Using (\ref{eqn:3-1}), we see for any integer $m$ that 
\[
\sum_{j=0}^{D} A_j (u_{m+j}+\epsilon_{m+j})=0
\]
and that 
\begin{align}\label{eqn:3-2}
s_m=-A_0 \epsilon_{m}-A_1\epsilon_{m+1}-\cdots-A_D\epsilon_{m+D}.
\end{align}
Hence, $(s_m)_{m\in \Z}$ is a bounded sequence of integers. In fact, 
putting $B:=(|A_0|+|A_1|+\cdots+|A_D|)/2$, we have 
$|s_m|\leq B$ for any $m\in \Z$. Note that if $m$ is sufficiently large, 
then $u_{-m}=0$, and hence, $s_{-m}=0$ for any sufficiently large $m\in \Z$. 
Thus, $(s_m)_{m\in \Z}\in\Omega_0$. 
\par
In what follows, we represent 
$(\epsilon_m)_{m\in\Z}$ by $(s_m)_{n\in\Z}$. 
We use the letter $\1$ to denote an indicator function. 
For instance, $\1_{m=n}$ is 1 when $m=n$ and $0$ otherwise.
We consider square matrices 
$\mathcal{A}=(a_{m,n})_{m,n\in \Z}, \mathcal{B}=(b_{m,n})_{m,n\in \Z}$ 
with the index set $\Z\times \Z$, 
where $m$ and $n$ denote 
the row number and the column number, respectively. 
When $\sum_{j\in \Z} a_{m,j} b_{j,n}$ converge for any $m,n\in \Z$, 
we define the product of $\mathcal{A}, \mathcal{B}$ by 
\begin{align}\label{eqn:3-4}
(a_{m,j})_{m,j\in \Z} (b_{j,n})_{m,n\in \Z}=
\left(\sum_{j\in \Z} a_{m,j} b_{j,n}\right)_{m,n\in \Z}.
\end{align}
We call that $\mathcal{B}$ an inverse matrix of $\mathcal{A}$ if 
both of $\mathcal{A}\mathcal{B}$ and $\mathcal{B}\mathcal{A}$ are 
defined and if \[
\mathcal{A}\mathcal{B}=\mathcal{B}\mathcal{A}
=(\1_{m=n})_{m, n\in \Z}. \]
We also consider the column vector $\by=(y_m)_{m\in \Z}$ 
with the index set $\Z$. 
When $\sum_{j\in \Z} a_{m,j} y_j$ converges for any $m\in \Z$, then 
we define the product $\mathcal{A}\by$ by 
\begin{align}\label{eqn:3-5}
(a_{m,j})_{m,j\in \Z} (y_j)_{j\in \Z}=
\left(
\sum_{j\in \Z} a_{m,j} y_j
\right)_{m\in \Z}.
\end{align}
Note that even if both of $(\mathcal{A} \mathcal{B}) \by$ and 
$\mathcal{A} (\mathcal{B} \by)$ are defined, these two vectors 
do not necessarily coincide with each other. \par
We consider the column vector $\bm{\epsilon}=(\epsilon_n)_{n\in \Z}\in\bE_f$. 
\textcolor{red}{Defining 
a column vector $\bs$ by $\bs=(s_m)_{n\in\Z}$, we see by (\ref{eqn:3-2}) that $\bs=\mathcal{A}_f\bep$,}
where the square matrix $\mathcal{A}_f$ with the index set 
$\Z\times \Z$ is defined as follows: 
\begin{align*}
\mathcal{A}_f:=(-A_{n-m} \1_{0\leq n-m\leq D})_{m,n\in \Z}.
\end{align*}
This is a matrix representation of the relation (\ref{eqn:3-2}):
\[
s_m=-\sum_{n=m}^{m+D}A_{n-m}\epsilon_n
=\sum_{n\in \Z}(-A_{n-m} \1_{0\leq n-m\leq D})\epsilon_n.
\]
We shall represent $\bep$, using the inverse matrix of $\mathcal{A}_f$. 
For any integer $n$, put 
\[
\rho^f_n=\frac{-1}{2\pi i}\int_{|z|=1+0}\frac{z^{n-1}}{f(z)}dz,
\]
where the integral is around the unit circle to the positive direction, 
avoiding poles, if exist, from outside. 
It is easily seen that $\rho^f_n$ is a real number for any $n\in \Z$. 
We now define the square matrix $\mathcal{B}_f$ with the index set 
$\Z\times \Z$ by 
\begin{equation}\label{eqn000}
\mathcal{B}_f=(\rho_{m-n}^{f})_{m,n\in \Z}.
\end{equation}
\begin{theorem}\label{thm:3-1}
(1) $\mathcal{B}_f$ is the inverse matrix of $\mathcal{A}_f$.
\\
(2) If $(s_m)$ is defined in (\ref{eqn:3-2}) with respect to $(\epsilon_m)$ 
(\ref{eqn:3-0}), then we have
\[(\epsilon_m)=\mathcal{B}_f (s_m).\]
(3) \textcolor{red}{Suppose that $P(X)$ is monic. Then for}
any $(t_m)\in\Omega_0$, there exists $(\epsilon_m)\in\bE_f$ such that 
$$\mathcal{B}_f (t_m)\equiv (\epsilon_m) \pmod{\Z}.$$
\end{theorem}

For the calculation of $\epsilon_m$, see also (\ref{eqn:3-17.5}). Before proving Theorem \ref{thm:3-1},
note that $\mathcal{A}_f\mathcal{B}_f$ and $\mathcal{B}_f\mathcal{A}_f$ 
are well-defined because their components are obtained as finite sums. 
We need the following lemmas. 

\begin{lem}\label{lem:3-1}
For any integer $n$, we have 
\[
\sum_{j=0}^D -A_j\rho^f_{j+n}=\1_{n=0}.
\]
\end{lem}
\begin{proof}
\begin{align*}
&\sum_{j=0}^D -A_j\rho^f_{j+n}
=\sum_{j=0}^D A_j\frac{1}{2\pi i}\int_{|z|=1+0}\frac{z^{j+n-1}}{f(z)}dz\\
&=\frac{1}{2\pi i}\int_{|z|=1+0}\frac{\sum_{j=0}^D A_jz^j}{f(z)}z^{n-1}dz\\
&=\frac{1}{2\pi i}\int_{|z|=1+0}z^{n-1}dz=\1_{n=0}.
\end{align*}
\end{proof}

\begin{lem}\label{lem:3-3}
If $P(X)$ is monic, then 
\[\Res \left(\frac{z^{n-1}}{f(z)},\infty\right)\in \Z\]
for any integer $n$. 
\end{lem}
\begin{proof}
We use the same notation as the proof of Proposition \ref{prop:3-1}. 
If $P(X)$ is monic, then $f(X)$ is also monic. Note that 
\[
\frac{1}{X^Df(1/X)}=\frac{1}{1+A_{D-1}X+\cdots+A_0 X^D}
\in \Z[[X]]
\]
as formal power series. 
Thus, combining (\ref{eqn:3-11}) and (\ref{eqn:3-12}), 
we verify Lemma \ref{lem:3-3}. 
\end{proof}

\begin{lem}\label{lem:3-2}
(1) 
There exists a real number $\lambda>1$ such that 
$\rho_{-n}^{f}=o(\lambda^{-n})$ as $n$ tends to infinity. \\
(2) $\rho_n^{f}=O(n^{k})$ as n tends to infinity. 
Moreover, if $P$ is hyperbolic, then there exists a positive real number 
$\delta<1$ such that $\rho_n^{f}=o(\delta^n)$ as $n$ tends to infinity. 
\end{lem}

We shall verify Lemma \ref{lem:3-2} in Subsection \ref{subsection3-2}. 
See Remark \ref{rem:3-1}. 

We now come back to the proof of Theorem \ref{thm:3-1}. 
We calculate the $(m,n)$-component of $\A_f\B_f$ 
using Lemma \ref{lem:3-1} as follows: 
\begin{align*}
&(\A_f\B_f)_{m,n}=\sum_{j\in \Z}-A_{j-m} \1_{0\leq j-m\leq D}\cdot \rho_{j-n}^f
=\sum_{h\in \Z} -A_h \1_{0\leq h\leq D}\cdot \rho_{h+m-n}^f\\
&=\sum_{h=0}^D -A_h \rho_{h+m-n}^f=\1_{m-n=0}.
\end{align*}
Thus, we proved $\mathcal{A}_f\mathcal{B}_f=(\1_{m=n})_{m,n\in \Z}$. 
In the same way as above, we can show that 
$\mathcal{B}_f\mathcal{A}_f=(\1_{m=n})_{m,n\in \Z}$, 
which implies (1) of Theorem \ref{thm:3-1}. 
Since $(s_m)=\mathcal{A}_f\cdot(\epsilon_m)$, to prove (2) of Theorem \ref{thm:3-1}, 
it suffices to verify 
\begin{align}\label{eqn:3-6}\mathcal{B}_f (\mathcal{A}\cdot(\epsilon_m))=
(\mathcal{B}_f \mathcal{A}_f)\cdot(\epsilon_m).\end{align}
In fact, we see 
\[(\mathcal{B}_f \mathcal{A}_f)\cdot(\epsilon_m)
=(\1_{m=n})_{m,n\in \Z}\cdot(\epsilon_m)=(\epsilon_m).\]
Let $\mathcal{A}_f=:(a_{m,n})_{m,n\in \Z}$ and 
$\mathcal{B}_f=:(b_{m,n})_{m,n\in \Z}$. For the proof of (\ref{eqn:3-6}), 
it suffices to check for any $m\in \Z$ that 
$\sum_{k,j\in \Z}|b_{m,k}a_{k,j} \epsilon_j|<\infty$. 
By Lemma \ref{lem:3-2}, we see for any integer $\ell$ that 
\[
\sum_{j\in \Z} |\rho_{\ell-j}^f \epsilon_j|<\infty.
\]
Putting $A:=\max\{|A_j|\mid 0\leq j\leq D\}$, we get 
\begin{align*}
\sum_{k,j\in \Z}|b_{m,k}a_{k,j} \epsilon_j|
&=
\sum_{k,j\in \Z}|\rho_{m-k}^f A_{j-k} \1_{0\leq j-k\leq D}\cdot \epsilon_j|\\
&\leq A \sum_{j\in \Z}\sum_{k\in \Z}
|\rho_{m-k}^f \1_{0\leq j-k\leq D}\cdot \epsilon_j|\\
&=A\sum_{j\in \Z}\sum_{h\in \Z}
|\rho_{h+m-j}^f \1_{0\leq h\leq D}\cdot \epsilon_j|\\
&=A\sum_{h=0}^D\sum_{j\in \Z}
|\rho_{h+m-j}^f \cdot \epsilon_j|<\infty, 
\end{align*}
which implies (2) of Theorem \ref{thm:3-1}

Let us prove (3) of Theorem \ref{thm:3-1}. By (3) of Proposition \ref{prop:3-1}, 
there exists $\bh\in\Xi_k$ such that 
$(\rho_n^f)_{n\in\Z}\equiv (\epsilon({x_n}(\bh)))_{n\in\Z}\pmod\Z$. 
Take an arbitrary $(t_n)_{n\in\Z}\in\Omega_0$ such that $t_n=0~(\forall n<n_0)$. 
Then, we have
\begin{align*}
&(\mathcal{B}_f\cdot(t_n)_{n\in\Z})_m=\sum_{n=n_0}^\infty\rho_{m-n}^ft_n
\equiv\sum_{n=n_0}^\infty\sum_{j=1}^pg_j(m-n)\alpha_j^{m-n}t_n\\
&=\sum_{j=1}^p\sum_{n=n_0}^\infty\alpha_j^{-n}t_ng_j(m-n)\alpha^m
=\sum_{j=1}^ph_j(m)\alpha_j^m\pmod\Z,
\end{align*}
where $(g_1,\cdots,g_p)\in\Sigma^k$ with $g_j(X)=\sum_{\ell=0}^kc_{j\ell}X^\ell$ and 
$$
h_j(X)=\sum_{i=0}^k\left(\sum_{\ell=i}^k\sum_{n=n_0}^\infty c_{j\ell}\binom{\ell}{i}
(-n)^{\ell-i}\alpha_j^{-n}t_n\right)X^i. 
$$
Since $\bh:=(h_1,\cdots,h_p)\in\Xi_k$, 
$\mathcal{B}_f(t_n)_{n\in\Z}\equiv\epsilon(x(\bh))\in\bE_f\pmod\Z$, 
which completes the proof of Theorem \ref{thm:3-1}.

\subsection{Calculation of $\rho_n^f$}\label{subsection3-2}
Let $r(z)\in \C(z)$ be a rational function. 
For any $\lambda\in \C\cup \{\infty\}$, we denote the residue 
of $r(z)$ at $z=\lambda$ by $\Res (r(z),\lambda)$. 
\begin{prop}\label{prop:3-1}
(1) For any $n\geq 1$, we have 
\begin{align}\label{eqn:3-7}
\rho_n^f=-\sum_{j=p+1}^{d} \Res \left(\frac{z^{n-1}}{f(z)},\alpha_j\right).
\end{align}
(2) For any $n\leq D-1$, we have 
\begin{align}\label{eqn:3-8}
\rho_n^f=\sum_{j=1}^{p} \Res \left(\frac{z^{n-1}}{f(z)},\alpha_j\right).
\end{align}
(3) If $P(X)$ is monic, then 
$$\rho_n^f\equiv\sum_{j=1}^{p} \Res \left(\frac{z^{n-1}}{f(z)},\alpha_j\right)
\pmod\Z~~(\forall n\in\Z).$$
(4) If $P(X)$ is monic, then $(\epsilon(\rho_n^f))_{n\in\Z}\in\bE_f$. 
\end{prop}
\begin{proof}
By the definition of $\rho_n^f$, we have 
\begin{align}\label{eqn:3-9}
\rho_n^f=-\sum_{j=p+1}^{d} \Res \left(\frac{z^{n-1}}{f(z)},\alpha_j\right)
-\Res \left(\frac{z^{n-1}}{f(z)},0\right).
\end{align}
If $n\geq 1$, then (\ref{eqn:3-9}) implies (\ref{eqn:3-7}) 
by the assumption $f(0)\ne 0$. 
Since the sum of the residues of any rational function 
including the infinity is 0, we see 
\begin{align}\label{eqn:3-10}
\rho_n^f=\sum_{j=1}^{p} \Res \left(\frac{z^{n-1}}{f(z)},\alpha_j\right)
+\Res \left(\frac{z^{n-1}}{f(z)},\infty\right).
\end{align}
Putting 
\begin{align}\label{eqn:3-11}
\widetilde{f_n}(w):=\frac{(1/w)^{n-1}}{f(1/w)}(-w^{-2})
=-\frac{w^{-n+D-1}}{w^Df(1/w)},
\end{align}
we have 
\begin{align}\label{eqn:3-12}
\Res \left(\frac{z^{n-1}}{f(z)},\infty\right)
=\Res (\widetilde{f_n}(w),0).
\end{align}
Hence, if $n\leq D-1$, then we get 
\[
\Res \left(\frac{z^{n-1}}{f(z)},\infty\right)
=0,
\]
which implies (\ref{eqn:3-8}) by (\ref{eqn:3-10}). 

Let $1\leq j\leq d$. Then there exists a polynomial $r_{k,j}(z)$ of 
degree $k$ such that 
\begin{align}\label{eqn:3-13}
\Res \left(\frac{z^{n-1}}{f(z)},\alpha_j\right)
=r_{k,j}(n)\alpha_j^n
\end{align}
for any integer $n$. In fact, putting 
\[
g(z):=\frac{(z-\alpha_j)^{k+1}}{f(z)},
\]
where $g(\alpha_j)\ne 0,\infty$, we see 
\begin{align*}
\Res \left(\frac{z^{n-1}}{f(z)},\alpha_j\right)
&=
\Res \left(
(z-\alpha_j)^{-k-1}z^{n-1} g(z),\alpha_j
\right)\\
&=
\left. \frac{1}{k!}\frac{d^k}{d z^k} (z^{n-1} g(z))\right|_{z=\alpha_j},
\end{align*}
which implies (\ref{eqn:3-13}). In the case of $k=0$, then 
(\ref{eqn:3-13}) implies 
\begin{align}\label{eqn:res_k_0}
\Res \left(\frac{z^{n-1}}{f(z)},\alpha_j\right)
=\frac{\alpha_j^{n-1}}{f'(\alpha_j)}.
\end{align}

(3) follows from (2), (\ref{eqn:3-10}) and Lemma \ref{lem:3-3}.

Let us prove (4). By (3), it is sufficient to prove that 
$$
\Res \left(\frac{z^{n-1}}{f(z)},\alpha_j\right)=g_j(n)\alpha_j^n~~
(n\in\Z,~j=1,\cdots,p). 
$$
with $(g_1,\cdots,g_p)\in\Xi_k$, which follows from (\ref{eqn:3-13}) 
together with 
$$ 
\Res \left(\frac{z^{n-1}}{f(z)},\overline{\alpha}_j\right)
=\overline{r}_{k,j}(n)\overline{\alpha}_j^n.
$$
\end{proof}

\begin{rem}\label{rem:3-1}
Combining Proposition \ref{prop:3-1} and (\ref{eqn:3-13}), we 
get Lemma \ref{lem:3-2}. 
\end{rem}

\section{Applications to Markoff-Lagrange spectrum } \label{section4}
Now we describe the main results of our paper, derived as 
applications of our intertwining formula.

Recall for any $\bg\in \Xi_k$, the linear reccrence $\bx(\bg)=(x_n(\bg))_{n\in\Z}$ is defined by (\ref{basicdef}). Our interest is the maximal limit point $\limsup_{n\to\infty}\|x_n(\bg)\|$. 
We set 
\[
\LL^{(k)}(P):=\left\{\left.
    \limsup_{n\to\infty}\|x_n(\bg)\|\ \right| \bg \in \Xi_k
\right\}.
\]
In the case of $k=0$, we put  $\LL(P):=\LL^{(0)}(P)$. 
\begin{dfn}
Let $P(X)\in \Z[X]$ be as above. \\
(1) We call $P(X)$ hyperbolic if $|\alpha_i|\ne 1$ for any 
$1\leq i\leq d$. \\
(2) We call $P(X)$ expansive if $p=d$. 
\end{dfn}

\begin{theorem}\label{thm:2-1}
Suppose that $P(X)\in \Z[X]$ is monic and hyperbolic. 
Then, for any $k\geq 0$, the set $\LL^{(k)}(P(X))$ is closed in $[0,1/2]$. 
\end{theorem}

Pursuing the analogy to the Markoff-Lagrange spectrum, we expect that the set $\LL^{(k)}(P(X))$ is ``sparse" at $0$ and ``thick" around $1/2$
for hyperbolic $P(X)$. We confirm this speculation in some special cases.

\begin{theorem}\label{thm:iso}
Let $k$ be a nonnegative integer. 
Suppose that $P(X)$ is monic. 
Then 0 is an isolated point of 
the set $\LL^{(k)}(P(X))$ if and only if $P(X)$ is hyperbolic. Moreover, if $P(X)$ is not hyperbolic, then $\LL^{(k)}(P(X))$ is dense in 
$[0,1/2]$. 
\end{theorem}

Let $\lceil x \rceil$
be the smallest integer not less than $x$. 

\begin{theorem}\label{thm:pos_int_gen}
Suppose that $P(X)\in \Z[X]$ is monic and hyperbolic. 
Let 
\begin{align*}
\beta_1:&=
\min\{|\alpha_1|,|\alpha_2|,\ldots,|\alpha_p|\}\\
\beta_2:&=\min\{|\alpha_{p+1}|^{-1},\ldots, 
|\alpha_d|^{-1}\}\\
\beta:&=\min\{\beta_1,\beta_2\}
\end{align*}
(1) $1/2$ is an accumulation point of $\LL^{(k)}(P(X))$.\\
(2) Let $k=0$. Assume that either $\beta_1=\beta$ is attained by $|\alpha_i|$ at a single index $i=1,\dots,p$, or $\beta_2=\beta$ is attained by $|\alpha_i|^{-1}$ at 
a single index $i=p+1,\dots,d$, and 
\begin{align}
\label{eqn:3-19}
\left\lceil\frac{\beta-1}{2}\right\rceil\sum_{n\in \Z}|\rho_n^f|<\frac 12
\end{align}
holds, then there exists a sequence of proper
intervals $[c_i,d_i]$ contained in $\LL^{(0)}(P(X))$ and $\lim_{i\to \infty}c_i=1/2$.\\
(3) Let $k=0$. Assume that $\beta_1$ is attained by $|\alpha_i|$ at a single index $i=1,\dots,p$ 
and that $\beta_2$ is attained by $|\alpha_i|^{-1}$ at 
a single index $i=p+1,\dots,d$. 
Moreover, suppose that 
\begin{align}
\label{eqn:3-newbeta}
\left\lceil\frac{\widetilde{\beta}-1}{2}\right\rceil\sum_{n\in \Z}|\rho_n^f|<\frac 12,
\end{align}
where $\widetilde{\beta}:=\max\{\beta_1,\beta_2\}$. 
Then there exists a $v<1/2$ such that $[v,1/2]\subset \LL^{(0)}(P(X))$.

\end{theorem}

We note that the second and last statements of Theorem \ref{thm:pos_int_gen} are applicable to obtain proper intervals contained in $\LL^{(k)}(P(X))$ for general $k$. 
In fact, if $k$ is a positive integer, then we have $\LL^{(0)}(P(X))\subset \LL^{(k)}(P(X))$. \par
The condition (\ref{eqn:3-19}) 
becomes difficult to check when $f$ has complex roots but there is an algorithm to check if it is valid. 
By introducing a threshold $n_0\ge 0$,
Proposition 
\ref{prop:3-1} implies  
\begin{align*}
\sum_{n\geq 1}|\rho_n^f|&\le
\sum_{1\le n\le n_0}
\left|\sum_{j=p+1}^{d} 
\Res \left(\frac{z^{n-1}}{f(z)},\alpha_j\right)\right|
+
\sum_{n>n_0}
\sum_{j=p+1}^{d} \frac{|\alpha_j^{n-1}|}{|f'(\alpha_j)|}\\
&\le
\sum_{1\le n \le n_0}
\left|\sum_{j=p+1}^{d} 
\frac{\alpha_j^{n-1}}{f'(\alpha_j)}\right|
+
\sum_{j=p+1}^{d} \frac{|\alpha_j|^{n_0}}{(1-|\alpha_j|)\prod_{i\neq j}|\alpha_j-\alpha_i|}
\end{align*}
and the similar inequality 
\begin{align*}
\sum_{n\leq 0}|\rho_n^f|&\le
\sum_{1-n_0\le n\le 0}
\left|\sum_{j=1}^{p} 
\frac{\alpha_j^{n-1}}{f'(\alpha_j)}\right|
+
\sum_{j=1}^{p} \frac{|\alpha_j|^{-n_0}}{(|\alpha_j|-1)\prod_{i\neq j}|\alpha_j-\alpha_i|}.
\end{align*}
Therefore taking $n_0$ sufficiently large, we can 
confirm (\ref{eqn:3-19}).
Letting $n_0=0$, we see
\begin{equation}
\label{Suff}
    \left\lceil\frac{\beta-1}2\right\rceil\sum_{j=1}^d\frac{1}{\big||\alpha_j|-1\big|\prod_{i\neq j}|\alpha_j-\alpha_i|}<\frac 12.
\end{equation}
implies (\ref{eqn:3-19}).

Let $\alpha>1$ be an algebraic integer of prime
degree $>2$. Then $\alpha^m$ with $m\in \N$
is an algebraic integer of the same degree as well.
We denote its minimal polynomial by $f^{(m)}$. Then 
$\beta$ (resp. $\alpha_j$) is replaced by $\beta^m$ (resp. $\alpha_j^m$), and we see that
there exists $m_0$ that $f^{(m)}$ satisfies (\ref{Suff})
for $m\ge m_0$. Thus there exists infinitely many 
polynomials which satisfies (\ref{eqn:3-19}).


Let 
\begin{align*}
E(X):=\frac{1-(1-X)\prod_{m=0}^{\infty}(1-X^{2^m})}{2X}
\end{align*}
and 
\begin{align*}
E^{(k)}(X):=\frac{1+X^{2^k}-(1-X)\prod_{m=0}^{k-1}(1-X^{2^m})}
{2X(1+X^{2^k})}
\end{align*}
for $k\geq 0$, where $\prod_{m=0}^{k-1}(1-X^{2^m})=1$ for $k=0$. 
\textcolor{blue}{
For any $m\geq 0$ and $r\geq 1$, we define the complete homogeneous symmetric polynomial by 
\begin{align*}
H_r^{(m)}(X_1,\ldots,X_r):=\sum_{
\begin{subarray}{c}
i_1,\ldots,i_r\geq 0\\
i_1+\cdots+i_r=m
\end{subarray}
}X_1^{i_1}\cdots X_r^{i_r}.
\end{align*}
It is known that 
\begin{align}\label{eqn:comppoly}
H_r^{(m)}(X_1,\ldots,X_r)
=
\sum_{i=1}^{r}\left(
\prod_{\begin{subarray}{c}1\leq j\leq r\\ j\ne i
\end{subarray}} \frac{X_i}{X_i-X_j}
\right)
X_i^m,
\end{align}
see \cite[Problem 7.4]{Stanley2}, \cite[Appendix A]{Louck-Biedenharn}, \cite[Theorem 1.7]{Gustafson-Milne} and \cite[Lemma 3.1]{kan2}.
Moreover, for any formal power series 
$F(X)=\sum_{m=0}^{\infty} a_m X^m\in \C[[X]]$, we have 
\begin{align}\label{eqn:homformal}
\sum_{m=0}^{\infty}a_m H_r^{(m)}(X_1,\ldots,X_r)
=
\sum_{i=1}^{r}\left(
\prod_{\begin{subarray}{c}1\leq j\leq r\\ j\ne i
\end{subarray}} \frac{X_i}{X_i-X_j}
\right)
F(X_i). 
\end{align}
In view of these formulas,} for any $r\geq 1$, we set 
\begin{align*}
E_r(X_1,\ldots,X_r)&:=\sum_{i=1}^{r}\left(
\prod_{\begin{subarray}{c}1\leq j\leq r\\ j\ne i
\end{subarray}} \frac{X_i}{X_i-X_j}
\right)
E(X_i)
,\\
E_r^{(k)}(X_1,\ldots,X_r)&:=\sum_{i=1}^{r}\left(
\prod_{\begin{subarray}{c}1\leq j\leq r\\ j\ne i
\end{subarray}} \frac{X_i}{X_i-X_j}
\right)
E^{(k)}(X_i).
\end{align*}
\begin{theorem}\label{thm:disc1}
Assume that $P(X)$ is monic and expanding. Moreover, suppose that 
$$\mbox{either }\rho_0^f\ge 2\rho_{-1}^f\ge 2^2\rho_{-2}^f\ge\cdots>0
\mbox{ or }\rho_0^f\le 2\rho_{-1}^f\le 2^2\rho_{-2}^f\le\cdots<0,$$
which is equivalent to that either
\begin{align}\label{eqn:assum_disc1}
&H_d^{(0)}(\alpha_1^{-1},\ldots,\alpha_d^{-1})
\ge 2H_d^{(1)}(\alpha_1^{-1},\ldots,\alpha_d^{-1})
\ge 2^2H_d^{(2)}(\alpha_1^{-1},\ldots,\alpha_d^{-1})
\ge\cdots>0\nonumber\\
&\mbox{or}\nonumber\\
&H_d^{(0)}(\alpha_1^{-1},\ldots,\alpha_d^{-1})
\le 2H_d^{(1)}(\alpha_1^{-1},\ldots,\alpha_d^{-1})
\le 2^2H_d^{(2)}(\alpha_1^{-1},\ldots,\alpha_d^{-1})
\le\cdots<0,
\end{align}
where $d$ is the degree of $P(X)$.
Then $|a_0|^{-1}E_d(\alpha_1^{-1},\ldots,\alpha_d^{-1})$ is the minimal 
limit point of $\LL(P(X))$. 
Moreover, we have 
\begin{align*}
&\LL(P(X))\cap \left(
0,\frac1{|a_0|}E_d(\alpha_1^{-1},\ldots,\alpha_d^{-1})
\right)\\
&\hspace{10mm} =
\left\{
\left.
\frac1{|a_0|}E_d^{(k)}(\alpha_1^{-1},\ldots,\alpha_d^{-1})
\ \right| \ 
k=0,1,\ldots
\right\}.
\end{align*}
\end{theorem}
\begin{cor}\label{cor:disc2}
Assume that $P(X)$ is monic and 
that $\alpha_j>1$ for any $1\leq j\leq d$. 
Moreover, suppose that 
\begin{align}\label{eqn:cor}
    \alpha_1^{-1}+\alpha_2^{-1}+\cdots+
\alpha_d^{-1}\leq \frac12. 
\end{align}
Then $|a_0|^{-1}E_d(\alpha_1^{-1},\ldots,\alpha_d^{-1})$ is the minimal 
limit point of $\LL(P(X))$. 
Moreover, we have 
\begin{align*}
&\LL(P(X))\cap \left(
0,\frac1{|a_0|}E_d(\alpha_1^{-1},\ldots,\alpha_d^{-1})
\right)\\
&\hspace{10mm} =
\left\{
\left.
\frac1{|a_0|}E_d^{(k)}(\alpha_1^{-1},\ldots,\alpha_d^{-1})
\ \right| \ 
k=0,1,\ldots
\right\}.
\end{align*}
\end{cor}
Corollary \ref{cor:disc2} follows from Theorem  \ref{thm:disc1}. In fact, Since 
$\alpha_1,\ldots,\alpha_d>1$,  
assumption (\ref{eqn:cor}) implies for any $m\geq 0$ that 
\begin{align*}
&H_d^{(m+1)}(\alpha_1^{-1},\ldots,\alpha_d^{-1})\\
&\leq 
(\alpha_1^{-1}+\alpha_2^{-1}+\cdots+
\alpha_d^{-1})
H_d^{(m)}(\alpha_1^{-1},\ldots,\alpha_d^{-1})\\
&\leq 
\frac{1}{2}
H_d^{(m)}(\alpha_1^{-1},\ldots,\alpha_d^{-1}).
\end{align*}

\begin{rem}
\textcolor{blue}{
Our definitions of $E(X)$ and $E^{(k)}(X)$ came from the study of
Dubickas \cite{dub2} on the simplest case:}
\[
\LL(X-a)=\left\{
\limsup_{n\to\infty}\|\xi a ^n\|\mid \xi\in \R
\right\},
\]
where $a$ is an integer greater than 1.
Introducing the function $E(X)$, 
Dubickas \cite{dub2} essentially proved that $a^{-1} E(a^{-1})$ is the minimal limit point of $\LL(X-a)$. \textcolor{blue}{From his proof, it is not 
difficult to describe 
all points in 
$\LL(X-a)\cap (0,a^{-1} E(a^{-1}))$ in terms of $E^{(k)}(X)$.
See the proofs of Lemma 4 in \cite{dub2} and Theorem 3.1 in \cite{Ak-Ka:21}.}
In \cite{Ak-Ka:21}, $E^{(k)}(X)$ was \textcolor{blue}{made explicit.}
For instance, we have 
\begin{align*}
E^{(0)}(X)=1-X+X^2-X^3+X^4-\cdots=\frac{1}{1+X}.
\end{align*}
\end{rem}

\begin{rem}
Let $\alpha>1$ be an algebraic number of 
degree $d\geq 2$. Let $\alpha_1=\alpha,\alpha_2,\ldots,\alpha_d$ be the conjugates of $\alpha$. Assume that $\alpha_1,\alpha_2,\ldots,\alpha_d>1$ and (\ref{eqn:cor}) holds. Using Theorem 2.3 in \cite{kan1}, we can show for any real $\xi\ne 0$ that 
\begin{align*}
    \limsup_{n\to\infty}\|\xi\alpha^n\|
    \geq \frac1{|a_0|}E_d(\alpha_1^{-1},\ldots,\alpha_d^{-1}),
\end{align*}
where $a_dX^d+\cdots+a_0\in \Z[X]$ is the minimal polynomial of $\alpha$. The conditions of Corollary 2.5 in \cite{kan1} should be corrected as above. 

\end{rem}

\section{Proof of Theorems in Section \ref{section4}}\label{section5}
\begin{proof}[Proof of Theorem \ref{thm:2-1}]
We can show Theorem \ref{thm:2-1} in the same way as the proof of 
the closedness of $\LL(\alpha)$ for a Pisot number $\alpha$. 
See the proof of Theorem 2.1 in \cite{Ak-Ka:21}. 
\end{proof}

\begin{proof}[Proof of Theorem \ref{thm:iso}]
First, we consider the case where $P(X)$ is hyperbolic. 
Suppose that $\bg\in \Xi_k$ satisfies 
\[
\limsup_{n\to\infty}\|x_n(\bg)\|
<\frac{1}{2|A_D|+2|A_{D-1}|+\cdots+2|A_0|}.
\]
Then (\ref{eqn:3-2}) implies that $|s_m|<1$ for any sufficiently large 
$m$ and that $s_m=0$ by $s_m\in \Z$. Recall that $s_{-m}=0$ for 
any sufficiently large $m$. Thus, 
Lemma \ref{lem:3-2} and the representation (\ref{eqn:3-17.5}) 
imply that $\limsup_{n\to\infty}\|x_n(\bg)\|=0$, which implies that 
$0$ is an isolated point of $\LL^{(k)}(P(X))$. \par
We now prove that 0 is an accumulation point of 
$\LL^{(k)}(P(X))$ in the case where $P(X)$ is not hyperbolic. 
We may assume that $k=0$ because 
$\LL^{(k)}(P(X))\supset \LL^{(0)}(P(X))$. 
Let $P(X)=P_1(X)P_2(X)\cdots P_{\ell}(X)$ be the factorization of 
$P(X)$ in $\Z[X]$. 
Recall for any $1\leq j\leq \ell$ that there exists at least one 
zero $\alpha$ of $P_{j}(X)$ with absolute value greater than 1. 
There exists $1\leq t\leq \ell$ such that $P_{t}(X)$ is not hyperbolic. 
Since $\LL(P(X))\supset \LL(P_{t}(X))$, 
we may assume that $P(X)=P_{t}(X)$. \par
When $\alpha$ is real, 
applying the Minkowski's convex body theorem, we can find a Pisot number $\eta$ with $\Q(\eta)=\Q(\alpha)$. 
In general, we can prove that there exists an algebraic integer $\eta$ satisfying: 
\begin{enumerate}
\item $|\eta|>1.$
\item Let $\eta'$ be a conjugate of $\eta$. Then $|\eta'|\geq 1$ if and only if $\eta'\in \{\eta,\overline{\eta}\}$. 
\item $\Q(\eta)=\Q(\alpha)$. 
\end{enumerate}
Let $\alpha_1=\alpha,\alpha_2,\ldots,\alpha_d$ 
be the Galois conjugates of $\alpha$ with 
$|\alpha_j|>1$ for $1\leq j\leq p$ and $|\alpha_h|\leq 1$ for 
$p+1\leq h\leq d$. 
Let $\psi^{(k)}$ be the embedding of 
$\Q(\alpha)$ to $\C$ with $\psi^{(k)}(\alpha)=\alpha_k$.
There exists an constant $0<\delta<1$ 
such that 
\begin{align*}
\limsup_{n\to\infty}\left\|
\sum_{j=1}^{p}\psi^{(j)}(\eta^M)\alpha_j^n
\right\|
=
\limsup_{n\to\infty}\left\|
\sum_{j=p+1}^{d}\psi^{(j)}(\eta^M)\alpha_j^n
\right\|
\leq d\delta^M.
\end{align*}
In what follows, we assume that $M$ satisfies 
$d\delta^{M}<1/2$. 
To check that $0$ is an accumulation point, 
it suffices to show that 
\[
\limsup_{n\to\infty}\left\|
\sum_{j=1}^{p}\psi^{(j)}(\eta^M)\alpha_j^n
\right\|>0.
\]
Suppose on the contrary that 
\[
\lim_{n\to\infty}\left\|
\sum_{j=1}^{p}\psi^{(j)}(\eta^M)\alpha_j^n
\right\|=0.
\]
Following discussion is found in Bugeaud \cite{Bu} in the proof of Hardy's theorem \cite{har}.
For simplicity, set $\lambda_j=\psi^{(j)}(\eta^M)$ 
for $j=1,2,\ldots,d$. 
Putting 
\[
u_n:=u\left(\sum_{j=1}^{p}\lambda_j\alpha_j^n\right), \quad
\epsilon_n:=\epsilon\left(\sum_{j=1}^{p}\lambda_j\alpha_j^n\right),
\]
we see 
\[
u_n=\sum_{j=1}^d\lambda_j\alpha_j^n 
, \quad
\epsilon_n=-\sum_{j=p+1}^d\lambda_j\alpha_j^n.
\]
We fix $k$ with $p+1\leq k\leq d$. Put 
\[
S(X)=
\sum_{m=0}^{d-1} c_m X^m
:=\prod_{1\leq j\leq d, j\ne k}(X-\alpha_j).
\]
Observe that 
\begin{align*}
& 
-\sum_{m=0}^{d-1} c_m \epsilon_{n+m}
=
\sum_{j=1}^{p}\lambda_j\alpha_j^n S(\alpha_j)
-\sum_{m=0}^{d-1} c_m \epsilon_{n+m}\\
& =
\sum_{m=0}^{d-1}c_m \sum_{j=1}^{p}   \lambda_j \alpha_j^{m+n}
-\sum_{m=0}^{d-1} c_m \epsilon_{n+m}\\
&=\sum_{m=0}^{d-1} c_m u_{n+m}
=\sum_{j=1}^d \lambda_j \alpha_j^{n} S(\alpha_j)=
\lambda_k \alpha_k^n S(\alpha_k).
\end{align*}
Hence, using $\lim_{n\to\infty} \epsilon_n=0$ and 
$\lambda_k S(\alpha_k)\ne 0$, we see that 
$|\alpha_k|<1$, which contradicts the assumption that 
$P(X)$ is hyperbolic because $k$ is any indices with $p+1\leq k\leq d$. 
\par
If $P(X)$ is not hyperbolic, then using the fact that $0$ is the accumulation point of $\LL^{(k)}(P(X))$, we can show that $\LL^{(k)}(P(X))$ is dense in $[0,1/2]$. 
In fact, take $\bg\in \Xi_k$ such that $\tau:=\limsup_{n\to\infty}\|x_n(\bg)\|$ is a sufficiently small positive real number. If a positive integer $h$ satisfies $h\tau<1/2$, then we see $h\tau=\limsup_{n\to\infty}\|x_n(h\bg)\|\in \LL^{(k)}(P(X))$. 
\end{proof}

\begin{proof}[Proof of Theorem \ref{thm:pos_int_gen}]
We show the first statement. 
Since $\LL^{(0)}(P(X))\subset \LL^{(k)}(P(X))$, 
we may assume that $k=0$, and so $f(X)=P(X)$. 
Let us assume that $P(X)$ is monic and hyperbolic. 
Defining the sequence $\bg_0\in \Xi$ by 
$$
\bg_0:=\begin{cases}
(-1/2,0,\ldots,0) & \mbox{ if }\alpha_1\in \R\\
(-1/4,-1/4,0,\ldots,0) & \mbox{ if }\alpha_1\not\in \R
\end{cases},
$$
we see $x_0(\bg_0)=-1/2$. For the convenience of notation, we define $\bt$ by $t_n=s_n(\bg_0)$ for any $n\in \Z$. Then Theorem \ref{thm:3-1} implies 
$$
\left| 
    \sum_{m\in \Z} \rho_{-m}^{f} t_m \right|
    =|\epsilon(x_0(\bg_0))|=1/2.
$$
Let $\delta$ be an arbitrary positive real number. 
Choosing a sequence $\bs=(s_n)_{n\in \Z}$ 
of the form 
\begin{equation}
    \label{PerApp}
\bs=0^{\infty}(t_{-R}t_{-R+1}\ldots t_0\ldots t_R
0^a 1 0^b)^{\infty},
\end{equation}
where $R$, $a$, and $b$ are suitable positive integers, we get 
\[
\frac12-\delta<\limsup_{n\to \infty}\left| 
    \sum_{m\in \Z} \rho_{n-m}^{f} s_m \right|
    <\frac12
\]
because $P(X)$ is hyperbolic. 
Consequently, $1/2$ is an accumulation point of 
$\LL^{(0)}(P(X))$, which implies the
first assertion. 

We now verify the second statement. 
By assumption, we have either 
\begin{align}\label{eqn:4-a-1}
\frac 1{\beta}=\lim_{n\to\infty}\left|\frac{\rho_{-n-1}^f}{\rho_{-n}^f}\right|
=\frac{1}{\min\{|\alpha_i|\mid 1\leq i\leq p\}}
\end{align}
at a single index $i=1,\ldots, p$, or 
\begin{align}\label{eqn:4-a-2}
\frac 1{\beta}=\lim_{n\to\infty}\left|
\frac{\rho_{n+1}^f}{\rho_n^f}
\right|
=
\frac{1}{\min\{|\alpha_i|^{-1}\mid p+1\leq i\leq d\}}
\end{align}
at a single index $i=p+1,\ldots,d$. 
We can easily show the following lemma: 
\begin{lem}\label{lem:4-a-1}
Let $(r_n)_n$ be a sequence of nonnegative real numbers such that 
$r_n\ne 0$ for any sufficiently large $n$. Suppose that 
\[
\lim_{n\to\infty}\frac{r_{n+1}}{r_n}\in (0,1).
\]
Take an integer $a$ greater than 1 with 
\begin{align*}
\frac{1}{a}<\lim_{n\to\infty}\frac{r_{n+1}}{r_n}.
\end{align*}
Let $\lfloor x\rfloor$ be the greatest integer not exceeding $x$
and $A:=\lfloor a/2\rfloor$. Then for any integer $N\geq 0$, the set 
\begin{align*}
\left\{\left.
\sum_{n=N}^{\infty} b_n r_n \ \right| \ 
b_n\in \Z, |b_n|\leq A\mbox{ for any }n\geq N
\right\}
\end{align*}
contains a proper interval containing $0$ as an inner point. 
\end{lem}

Put $a:=1+\lfloor \beta\rfloor(\geq 2)$ and 
$A:=\lfloor a/2\rfloor$. We prove 
the second assertion of Theorem \ref{thm:pos_int_gen} for the case (\ref{eqn:4-a-1}).
Lemma \ref{lem:4-a-1} implies for any integer $M,N\geq 0$ with 
$M<N$ that the set 
\begin{align}\label{eqn:set_m_n}
\mathcal{S}_{\pm}(A;M,N):=
\left\{\left.
\pm\rho_{-M}^f+
\sum_{n=N}^{\infty} b_n \rho_{-n}^f \ \right| \ 
b_n\in \Z, |b_n|\leq A\mbox{ for any }n\geq M
\right\}
\end{align}
contains a proper interval $J_{\pm}(M,N)$ with 
$\pm\rho_{-M}^f\in J_{\pm}(M,N)$. 
Suppose that $M$ is a sufficiently large integer. 
We observe that if 
$N$ is sufficiently large depending on $M$, then 
\begin{align}\label{eqn:4-a-3}
\left| \pm\rho_{-M}^f+\sum_{n=N}^{\infty} b_n \rho_{-n}^f \right|>
\left| \pm\rho_{-h-M}^f+\sum_{n=N}^{\infty} b_n \rho_{-h-n}^f \right|
\end{align}
for any positive integer $h$. 
Using a similar construction to 
(\ref{PerApp}), we can find an increasing sequence $(R_i)_i$ of integers and 
$\bu_{i}=(u_{i,n})_{n\in \Z}$ of the form
\begin{equation}
\label{Central}
\bu_{i}=
0^{\infty}t_{-R_i}t_{-R_i+1}\ldots t_0\ldots t_{R_i-1}t_{R_i} 0^{a_i}1 0^{\infty}
\end{equation}
with
$$
u_{i,0}u_{i,1}\dots u_{i,K_i}
=t_{-R_i}t_{-R_i+1}\ldots t_0\ldots t_{R_i-1}t_{R_i} 0^{a_i}1,
$$ 
where $K_{i}:=2R_i+a_i+1$. Noting for each $i$ that $u_{i,m}=0$ for all but finitely many $m$, we observe that 
$$
\max \left\{\left.\left\Vert\sum_{m\in \Z}\rho^f_{n-m}u_{i,m}\right\Vert \right| n\in \Z \right\}<\frac 12
$$
and
$$\lim_{i\to \infty}
\max \left\{\left.\left\Vert\sum_{m\in \Z}\rho^f_{n-m}u_{i,m}\right\Vert \right| n\in \Z \right\} =\frac 12.
$$
Let $i$ be any positive integer. There exists an integer $\ell$ with 
\begin{align}\label{eqn:4-a-4}
\left\Vert\sum_{m\in \Z}\rho_{\ell-m}^f u_{i,m}\right\Vert=\max\left\{
\left.\left\Vert\sum_{m\in \Z}\rho_{n-m}^f u_{i,m}\right\Vert\ \right|\ n\in \Z\right\}.
\end{align}
Let $\ell_{i}$ be the largest
$\ell$ which satisfies this equality. 
Putting 
$$
\frac 12-\delta_i:=
\left\Vert\sum_{m\in \Z}\rho_{\ell_{i}-m}^f u_{i,m}\right\Vert,
$$
we have $\delta_i>0$ and 
$\lim_{i\to \infty}\delta_i=0$ because $(\bu_i)_i$ approximates $\bt=(t_n)_{n\in \Z}$.


Let $M_{i}$ be a sufficiently large integer depending on $\ell_i$. 
We choose a sufficiently large integer $N_i$
depending $\ell_i$ and $M_i$.
Assuming (\ref{eqn:3-19}), we shall prove that 
\begin{align}\label{eqn:4-a-6}
\mathcal{I}_i:=\left(\frac 12-\delta_i+
J_{\pm}(M_{i},N_{i})\right)\cap\left[0,\frac12 \right]\subset\LL^{(0)}(P(X)).
\end{align}
\begin{rem}\label{rem:interval}
If there exist an index $i$ and a real number $v$ with $0<v<1/2$ such that 
\begin{align}\label{eqn:intervalhalf}
    \left[v,\frac12 \right]
    \subset\frac12 -\delta_i+J_{\pm}(M_{i},N_{i}),
\end{align}
then we have $[v,1/2]\subset \LL^{(0)}(P(X))$. 
In other words, an interval $c_i,d_i$ in the second statement satisfies $d_i=1/2$. 
We shall use the fact above to prove the last statement. 
\end{rem}
Recall that $J_{\pm}(M_{i},N_{i})$ is an interval contained in the set (\ref{eqn:set_m_n}).  
Let $z$ be any element of $J_{\pm}(M_i,N_i)$ with $1/2-\delta_i+z\in [0,1/2]$. Then $z$ can be represented as 
\[
z=\pm\rho_{-M_i}^f+\sum_{n\geq N_i}z_n\rho_{-n}^f
\]
with $|z_n|\le A$. 
For any integer $n$ with $n\geq N_i$, we define the finite word 
$v(i,n)$ by 
\[
v(i,n)=u_{i,0}u_{i,1}\dots u_{i,K_i}0^{M_i+\ell_i-K_i-1}
(\pm 1) 0^{N_i-M_i-1} z_{N_i}z_{N_i+1}\cdots z_n
\]
where we choose the signature of $\pm 1$ 
by the signature of $J_{\pm}(M_{i},N_{i})$.
We construct a sequence $\bw_i=(w_{i,n})_{n\in\Z}$ by 
$w_{i,n}=0$ for any $n\leq -1$ and 
\[
w_{i,0}w_{i,1}w_{i,2}\ldots=v(i,\kappa(1))0^{\tau(1)}v(i,\kappa(2))0^{\tau(2)}\cdots,
\]
where $(\kappa(m))_{m\geq 1}$ and $(\tau(m))_{m\geq 1}$ are 
sufficiently rapidly increasing sequence of positive integers. 
We take a suitable $\bg_i\in \Xi_0$ satisfying (\ref{eqn:3-18}). 
Recall that $\ell_i$ is the largest
integer satisfying (\ref{eqn:4-a-4}). 
Since $M_i, N_i$ are sufficiently large and 
$(\kappa(m))_{m\geq 1}$, $(\tau(m))_{m\geq 1}$ increase 
sufficiently rapidly, we obtain 
\begin{align*}
\limsup_{n\to\infty}\|x_n(\bg_i)\|
&=\limsup_{n\to\infty}|(\sigma^n(\bw_i))_f|\\
&=\left|\sum_{m\in \Z}\rho^f_{\ell_i-m}u_{i,m}+z\right|
\in \LL^{(0)}(P(X)),
\end{align*}
where $\sigma$ denotes the shift operator. 
Since $z$ is an arbitrary element of $J_{\pm}(M_i,N_i)$ with $1/2-\delta_i+z\in [0,1/2]$, we deduce (\ref{eqn:4-a-6}). 
We can show the case 
(\ref{eqn:4-a-2}) in the same way.


Now we prove the last statement. First we consider the case of (\ref{eqn:4-a-1}). 
By (\ref{eqn:4-a-1}) and (\ref{eqn:3-newbeta}), we can find a positive integer $m$ and an expansion
$$
\frac 12 =\sum_{n=1}^{\infty}
t_{-n}\rho_n^{f},
$$
where $(t_{-n})_{n\geq 1}$ is a sequence of integers such that 
$|t_{-k}|\le \lfloor (\widetilde{\beta}-1)/2 \rfloor$
for any $k\ge m$.
Instead of (\ref{Central}), we can take
$$
\bu_i=0^{\infty}10^{a_i}t_{-R_i}t_{-R_i+1}\ldots t_0 0^{\infty}
$$
using different indices:
\textcolor{red}{
$$
u_{i,-K_i}u_{i,1-K_i}\ldots u_{i,0}
=10^{a_i}t_{-R_i}t_{-R_i+1}\ldots t_0
$$
}
and $K_i=2+a_i+R_i$. 
Let 
$$x:= \left\lceil\frac{\widetilde{\beta}-1}{2}\right\rceil\sum_{n\in \Z}|\rho_n^f|, \quad y:=\frac12-x(>0).$$
Defining $\gamma_i$ by 
$$ \gamma_i:=\max\left\{
\left.\left\Vert\sum_{m\in \Z}\rho_{n-m}^f u_{i,m}\right\Vert\ \right|\ n\in \Z\right\}, $$
we have $\gamma_i<1/2$ and $\lim_{i\to\infty}\gamma_i=1/2$. Note that $x+y/2<1/2$. In what follows, we take a sufficiently large $i$ so that $\gamma_i>x+y/2$. 
Since $|t_{-k}|\le \lfloor (\beta-1)/2 \rfloor$ for any $k\ge m$, there exists a sufficiently large integer $L$ independent of $i$ such that, for any $n$ with $|n|\geq L$, we have 
$$ \left\Vert\sum_{m\in \Z}\rho_{n-m}^f u_{i,m}\right\Vert <x+\frac12 y(<\gamma_i). $$
Hence, we can choose a bounded sequence $(\ell_i)_i$
which satisfies (\ref{eqn:4-a-4}). \par
In the similar way as the proof of the second statement, we can construct a sequence of proper intervals $[c_i,d_i]\subset \LL^{(0)}(P(X))$. 
Let $\widetilde{A}:=
\lfloor(1+\lfloor\widetilde{\beta}\rfloor)/2\rfloor$. 
Let $J_{\pm}(M_i,N_i)$ be the interval in 
$\mathcal{S}_{\pm}(\widetilde{A};M_i,N_i)$, which is defined by (\ref{eqn:set_m_n}). 
Then $\mathcal{I}_i$ in (\ref{eqn:4-a-6}) gives a proper interval. 
Since $(\ell_i)_i$ is bounded, there exists a lower bound $z_0>0$ independent of $i$ satifying the following: $[c_i,d_i]$ fulfills $d_i-c_i>z_0$ unless $d_i=1/2$. 
We can construct an interval $[c_i,d_i]$ so that $c_i>1/2-z_0$, using the fact that, for any $N\geq 0$, we have $0$ is an inner point of the set 
\begin{align*}
    \left\{
    \left. \sum_{n=N}^{\infty} b_n\rho_{-n}^f
    \ \right| \ b_n\in\Z, |b_n|\leq \widetilde{A}
    \mbox{ for any }n\geq N
    \right\}
\end{align*}
by Lemma \ref{lem:4-a-1}. 
Therefore, we deduce that $d_i=1/2$, which implies the last statement. 
The case (\ref{eqn:4-a-2}) is
similar.
\end{proof}

\begin{rem}
\label{Interval}
The assumptions of Theorem \ref{thm:pos_int_gen}
are used to show that there exist positive integers $A$ and $m$,
so that
for any $|t_n|\le A$ with $t_n\in \Z$, we have
\begin{equation}
\label{SizeA}
A \sum_{n\in \Z} |\rho_{-n}^{f}|<\frac 12
\end{equation}
and
\begin{equation}
\label{NegInt}    
\left\{\left.
\sum_{n=m}^{\infty} t_{-n} \rho_{n}^{f}\ \right|
|t_n|\le A
\right\}
\end{equation}
contains $0$ as an inner point
and/or
\begin{equation}
\label{PosInt}
\left\{\left.
\sum_{n=m}^{\infty} t_n \rho_{-n}^{f}\ \right|
|t_n|\le A
\right\} 
\end{equation}
contains $0$ as an inner point. In cases,
we can show these assumptions by using the geometric property of the
attractor of iterated function system. 
For the first example in the introduction, we can choose $A=3$ to have (\ref{SizeA}) and
(\ref{NegInt}). For (\ref{PosInt}), we can confirm that 
the attractor $X$ of the iterated function system in the complex plane $\C$:
$$
X = \bigcup_{i=-3}^3 \frac1{\alpha} (X+i)
$$
is connected\footnote{One can also show that
$X$ is not simply connected and $0$ is not an inner point of $X$.}. This implies the projection of $X$ to the real axis is a proper interval containing the origin as an inner point. Thus
there exist an integer sequence $(t_n)_{n\ge 1}$ and
an integer $m$ that
$$
\frac 12=\sum_{n=1}^{\infty} t_n \rho_{-n}^f
$$
with $|t_n|\le A$ for $n\ge m$. Following the proof of Theorem \ref{thm:pos_int_gen}, 
there exists $v<1/2$ that $[v,1/2]\subset \LL_1$.
The idea to use fractal geometry may be extended in several directions but we do not explore more in this paper.
\end{rem}

\begin{proof}[Proof of Theorem \ref{thm:disc1}]
We can show Theorem \ref{thm:disc1} in the same method 
as the proof of Theorem 3.1 in \cite{Ak-Ka:21}. 
We give the sketch of the proof. 
By the symmetry, we may assume the first case of (\ref{eqn:assum_disc1}). 
In the proof, we put 
\[
e_k:=\frac{1}{|a_0|}E_d^{(k)}(\alpha_1^{-1},\ldots,\alpha_d^{-1})
\]
for $k\geq 0$
and 
\[
e:=\frac{1}{|a_0|}E_d(\alpha_1^{-1},\ldots,\alpha_d^{-1}).
\]
We make explicit our formula in the expanding case.
We have $\rho_n^f=0$ for any $n\geq 1$ by (\ref{eqn:3-7}). 
Moreover, if $n\leq 0$, then (\ref{eqn:3-8}) and (\ref{eqn:res_k_0}) 
imply that 
\begin{align}
\frac{1}{a_d}\rho_n^f
&=
\frac{1}{a_d}\sum_{i=1}^d\Res\left(\frac{z^{n-1}}{f(z)},\alpha_i\right)
=
\frac{1}{a_d}\sum_{i=1}^d\frac{\alpha_i^{n-1}}{f'(\alpha_i)}\nonumber\\
&=
\frac{1}{a_d}\sum_{i=1}^d\left(
\prod_{\begin{subarray}{c}1\leq j\leq d\\ j\ne i
\end{subarray}}
\frac{1}{\alpha_i-\alpha_j}\right)\alpha_i^{n-1}\nonumber\\
&=
\frac{1}{a_d\alpha_1\cdots\alpha_d}\sum_{i=1}^d
\left(
\prod_{\begin{subarray}{c}1\leq j\leq d\\ j\ne i
\end{subarray}}
\frac{1}{\alpha_j^{-1}-\alpha_i^{-1}}\right)\alpha_i^{n-d+1}\nonumber\\
&=
-\frac{1}{a_0}\sum_{i=1}^d
\left(
\prod_{\begin{subarray}{c}1\leq j\leq d\\ j\ne i
\end{subarray}}
\frac{\alpha_i^{-1}}{\alpha_i^{-1}-\alpha_j^{-1}}\right)
(\alpha_i^{-1})^{-n}
=-\frac{1}{a_0}H_d^{(-n)}(\alpha_1^{-1},\ldots,\alpha_d^{-1}).
\label{eqn:exp_explicit}
\end{align}
Let $\bg\in \Xi(=\Xi_0)$ and 
$s_n:=s_n(\bg)$. 
Then the second statement of 
Theorem \ref{thm:3-1} and (\ref{eqn:exp_explicit}) imply that
\begin{align}
\epsilon_n&=\sum_{j\in \Z} \rho_{-j}^f s_{j+n}
=\sum_{j=0}^{\infty}\rho_{-j}^f s_{j+n}\nonumber\\
&=-\frac{1}{a_0}\sum_{j=0}^{\infty} 
H_d^{(j)}(\alpha_1^{-1},\ldots, \alpha_d^{-1}) s_{j+n}\nonumber\\
&=:\sum_{j=0}^{\infty}\mu_j s_{j+n}=(s_{n}s_{n+1}
\ldots)_{\mu},
\label{eqn:4-b-1}
\end{align}
where $\mu_j=-a_0^{-1}H_d^{(j)}(\alpha_1^{-1},\ldots, \alpha_d^{-1})$ and 
$(\bt)_{\mu}=\sum_{j=0}^{\infty} \mu_j t_j$ for any infinite sequence 
$\bt=(t_n)_{n\geq 0}$. Note that (\ref{eqn:assum_disc1}) implies that 
\begin{align}\label{eqn:4-b-2}
0<\mu_{j+1}
\leq 
\frac{1}{2}
\mu_j
\end{align}
for any $j\geq 0$. 
We now represent limit points in Theorem \ref{thm:disc1}, 
using the representation $(\bt)_{\mu}$. 
Let $\iota$ be the substitution on the alphabet $\{0,1\}$ defined by 
$\iota(0)=1, \iota(1)=100$. We define the finite words $A_n$ 
($n=0,1,\ldots$) by $A_0=1$ and $A_n=\iota(A_{n-1})$ for $n\geq 1$. 
Let $\omega$ be the infinite sequence defined by 
\[\omega=\lim_{n\to\infty} A_n=10011100100\ldots.\] 
Moreover, let $A_n^{\infty}:=A_nA_n\ldots$ be the 
periodic sequence with period $A_n$. 
In the rest of the proof, we denote $\ov{b}:=-b$ for any integer $b$. 
Let $\Phi: \{0,1\}^{\infty}\to\{0,1,\ov{1}\}$ be a map defined as follows: 
\begin{align*}
\Phi(t_0t_1t_2t_3\ldots):=10^{t_0}\ov{1}0^{t_1}10^{t_2}\ov{1}0^{t_3}1\ldots.
\end{align*}
Then we see 
\begin{align}\label{eqn:4-b-3}
e_k
=(\Phi(A_k^{\infty}))_{\mu}
\end{align}
for any $k\geq 0$ and 
\begin{align}\label{eqn:4-b-4}
e
=(\Phi(\omega))_{\mu}=
(10\ov{1}1\ov{1}010\ov{1}01\ov{1}1\ldots)_{\mu}. 
\end{align}
In fact, if $d=1$, then (\ref{eqn:4-b-3}) and (\ref{eqn:4-b-4}) were 
proved in the paper \cite{Ak-Ka:21} (see the proof of Theorem 3.1). 
The case of general $d$ is deduced by 
(\ref{eqn:homformal}). \par
We take $\bg\in \Xi$ with 
\begin{align}\label{eqn:4-b-5}
\limsup_{n\to\infty}\|x_n(\bg)\|\in
\left(
0,e
\right].
\end{align}
Put 
\[
S:=\limsup_{n\to\infty}|s_n|.
\]
If $S=0$, then we have 
$\limsup_{n\to\infty}|s_n|=0$ by (\ref{eqn:4-b-1}), 
which contradicts (\ref{eqn:4-b-5}). \par
We prove that $S=1$. Assume on the contrary that $S\geq 2$. 
$S$ or $\ov{S}$ occurs infinitely many times in 
the sequence $(s_n)_{n\geq 0}=s_0s_1\ldots$. 
We assume for some integer $h$ with $-S+2\leq h\leq S$ that 
the finite word $S h$ occurs infinitely often in 
$(s_n)_{n\geq 0}$. 
Using (\ref{eqn:4-b-2}) and (\ref{eqn:4-b-4}), we get 
\[
\limsup_{n\to\infty}\|x_n(\bg)\|\geq 
(S \ \ov{S-2} \ \ov{S}^{\infty})_{\mu}\geq (1 0^{\infty})_{\mu}
>e,
\]
which contradicts (\ref{eqn:4-b-5}). 
Thus, each of the words 
$S h$ ($-S+2\leq h\leq S$) occurs at most finitely many times 
in $(s_n)_{n\geq 0}$. Similarly, each of the words 
$\ov{S} h$ ($-S\leq h\leq S-2$) appears at most finitely many times 
in $(s_n)_{n\geq 0}$. \par
Next, we assume that the word $S \ \ov{S-1}$ appears infinitely many 
times in $(s_n)_{n\geq 0}$. For any sufficiently large $n$, we have 
\[
\mu_{j} s_{n+j}+\mu_{j+1}s_{n+j+1}\geq -S\mu_j-(S-2)\mu_{j+1}.
\]
Hence, using (\ref{eqn:4-b-2}) and (\ref{eqn:4-b-4}), we obtain
\begin{align*}
\limsup_{n\to\infty}\|x_n(\bg)\|
&\geq
\left(S\ov{S-1}\left(\ov{S} \ \ov{S-2}\right)^{\infty}\right)_{\mu}
\\
&\geq (10\ov{1}1\ov{1}010^{\infty})_{\mu}>e,
\end{align*}
a contradiction. 
Thus, $S \ \ov{S-1}$ occurs at most finitely many times 
in $(s_n)_{n\geq 0}$. Similarly, $\ov{S} (S-1)$
appears at most finitely many times 
in $(s_n)_{n\geq 0}$. \par
Therefore, we see that $s_0s_1\ldots$ has the form 
$v(S\ov{S})^{\infty}$ with certain finite word $v$. 
Consequently, 
\[
\limsup_{n\to\infty}\|x_n(\bg)\|=S\cdot ((1\ov{1})^{\infty})_{\mu}>e,
\]
a contradiction. \par
Finally, we deduce that $S=1$. In the similar method as above, we 
can show that each of the following words occurs at most finitely 
many times in $(s_n)_{n\geq 0}$: 
$11,\ov{1}\ \ov{1}, 101,\ov{1}0\ov{1}, 00$. 
If necessary, considering $\sigma^{R}((s_n)_{n\geq 0})$ instead of 
$(s_n)_{n\geq 0}$ with suitable $R\geq 0$, we may assume that 
$(s_n)_{n\geq 0}$ has the form 
\[
10^{t_0}\ov{1}0^{t_1}10^{t_2}\ov{1}0^{t_3}1\ldots=\Phi((t_n)_{n\geq 0})
\]
with certain $(t_n)_{n\geq 0}\in \{0,1\}^{\infty}$. 
\par
Let $\bt=t_0t_1\ldots,\bt'=t_0't_1'\ldots$ be finite or infinite 
words on the alphabet $\{0,1\}$. Suppose that there exists $n\geq 0$ 
with $t_n\ne t_n'$. We define $\bt>\bt'$ if and only if 
\[
(-1)^h (t_h-t_h')>0,
\]
where $h=\min \{n\geq 0\mid t_n\ne t_n'\}$. Using (\ref{eqn:4-b-2}), 
we see that if $\bt,\bt'$ are infinite words on $\{0,1\}$ with 
$\bt>\bt'$, then 
\[
(\Phi(\bt))_{\mu}> (\Phi(\bt'))_{\mu}.
\]
Let $\mathcal{W}$ be the set of infinite words $\bt$ on $\{0,1\}$ 
satisfying the following: For any nonempty finite word $v$ on $\{0,1\}$ 
with $v>\omega$, the word $v$ occurs at most finitely many times 
in $\bt$. Then the set $\LL(P(X))\cap (0,e]$ is represented as 
\[
\left\{
\left.
\left(
\Phi\left(\limsup_{N\to\infty} \sigma^N(\bt)\right)
\right)_{\mu}
\ \right| \ 
\bt\in \mathcal{W}
\right\},
\]
where $\limsup_{N\to\infty} \sigma^N(\bt)$ is defined with respect 
to the order $>$. The following relation was proved in \cite{Ak-Ka:21}: 
\begin{align*}
\left\{
\left.
\Phi\left(\limsup_{N\to\infty} \sigma^N(\bt)\right)
\ \right| \ 
\bt\in \mathcal{W}
\right\}
=\{A_k^{\infty}\mid k=0,1,\ldots\}\cup \{\omega\}.
\end{align*}
Therefore, Theorem \ref{thm:disc1} follows from 
(\ref{eqn:4-b-3}) and (\ref{eqn:4-b-4}). 
\end{proof}
\section{Open problems}\label{section6}
We list several open problems related to the
multiplicative Lagrange spectrum. 
\begin{enumerate}
    \item Suppose that $P(X)$ is monic and not hyperbolic. Theorem \ref{thm:iso} implies for any $k\geq 0$ that $\LL^{(k)}(P(X))$ is dense in $[0,1/2]$. Does $\LL^{(k)}(P(X))$ contain a proper interval of $[0,1/2]$ ?
    \item Is $\LL^{(k)}(P(X))$ a closed set in the case where $P(X)$ is not monic? For instance, is $\LL^{(0)}(2X-3)=\{\limsup_{n\to\infty}\|\xi(3/2)^n\|\mid \xi\in \R\}$ a closed set? This problem is related to numerical systems proposed in \cite{Ak-Fr-Sa:08}.
    \item Let $\alpha>1$ be an algebraic integer with minimal polynomial $P(X)$. Then Theorem \ref{thm:2-1} implies for any $k\geq 0$ that $\LL^{(k)}(P(X))$ is a closed set. On the other hand, is  $\{\limsup_{n\to\infty}\|\xi\alpha^n\|\mid \xi\in \R\}$ a closed set? Moreover, it is also an interesting problem to find a nonzero polynomial $A(X)\in \R[X]$ such that $\limsup_{n\to\infty}\|A(n)\alpha^n\|$ is small (see also \cite{dub1}). 
    \item Do there exist a polynomial $P(X)$ and an integer $k\geq 0$ such that $\LL^{(k)}(P(X))$ contains an interval of positive length, but $[v,1/2]\not\subset \LL^{(k)}(P(X))$ for any real number $v$ with $0<v<1/2$?
    \item Investigate the ``discrete part" of $\LL^{(k)}(P(X))$ in the case where $P(X)$ is monic and hyperbolic. For instance, let $x$ be the minimal limit point of $\LL^{(k)}(P(X))$. Then is $x$ a left limit point? In other words, is there a strictly increasing sequence $(x_n)_{n\geq 1}$ of elements in $\LL^{(k)}(P(X))$ such that $\lim_{n\to\infty}x_n=x$?
    \item For any $S\subset \R$, we denote its Hausdorff dimension by $\dim_H(S)$. Let  $G(t):=\dim_H(\LL^{(k)}(P)\cap[0,t])$. Can we estimate $G(t)$ in the case where $P(X)$ is monic and hyperbolic? Moreover, is $G(t)$ continuous in $t$? In \cite{Moreira.18}, Moreira proved the continuity of the Hausdorff dimension for the classical Lagrange spectrum. 
\end{enumerate}

\section*{Acknowledgments}

This research was partially supported by JSPS grants (17K05159, 20K03528, 19K03439, 17H02849, 21H00989).

\end{document}